\documentclass[12pt]{article}
\pdfoutput=1
\usepackage{graphicx}
\usepackage{multicol}
\usepackage{wrapfig}

\newcommand{\beq}{\begin{eqnarray}}
\newcommand{\eeq}{\end{eqnarray}}
\newcommand{\nbeq}{\begin{eqnarray*}}
\newcommand{\neeq}{\end{eqnarray*}}

\newcommand{\be}[1]{\begin{equation} \label{#1}}
\newcommand{\ee}{\end{equation}}

\def\th{\theta}

\begin{document}

\begin{center}
{\Large\bf Monotone Empirical Bayes Estimators for the Reproduction Number in Borel-Tanner Distribution}
\end{center}
\begin{center}
{\sc George P. Yanev, Roberto Colson}\\
{\it The University of Texas Rio Grande Valley\\
Edinburg, Texas, USA\\}
e-mail: {\tt george.yanev@utrgv.edu}
\end{center}

\begin{abstract}
 We  construct a monotone version of an empirical Bayes estimator for the parameter of the Borel-Tanner distribution. Some properties of the estimator's regret risk are illustrated through
simulations.
\end{abstract}

\section{Introduction}

 The probability mass function (p.m.f.) of Borel-Tanner (BT) distribution is 
\begin{equation} \label{bt2}
p_r(x;\th)
     = c_r(x) \theta^{x-r}e^{-\theta x} \qquad x=r, r+1,\ldots,
\end{equation}
where $0<\theta<1$, $r$ is a positive integer, and $ c_r(x):=rx^{x-r-1}/(x-r)! $

The BT distribution  arises, for example, branching processes models and queueing theory. Originally (\ref{bt2}) was derived 
as the distribution of the number of customers served in a busy period of a single-server queuing process, started with $r$ customers and having traffic intensity $\th$, assuming Poisson arrivals and constant service time. Later BT distribution appeared in the theory of branching processes. If the number of offspring that an individual has is Poisson-distributed with offspring mean $0<\th<1$, then the total progeny of a Galton-Watson process starting with $r$ ancestors is a random variable with p.m.f. (\ref{bt2}). More recently, the distribution has been used to model a variety of real-world phenomena including: coalescence models (Aldous (1999)),  highway traffic flows (Koorey (2007)), propagation of internet viruses (Sellke et al. (2005)), cascading failures of energy systems (Ren et al. (2013)) and  herd size in finance modeling (Nirey et al. (2012)). 
Our interest in estimating $\theta$ stems from its role as the reproduction number of an epidemic infection modeled by a branching process (Farrington et al. (2003)).

In the context of branching processes, the parametric Bayesian statistical approach was first explored by Dion (1972) and  Jagers (1975), Section 2.13 (see also Guttorp (1991), Chapter 4). 
Adopting the Bayesian framework, suppose $\th \in \Omega$ is a realization of a random variable $\Theta$,
having a prior distribution $G$. It is well-known that, under the squared error loss, the value  $\th_G(x)$ of the Bayesian estimator for $\th$ is the posterior mean
\begin{equation} \label{bayes}
\th_G(x)  =   E\left[ \Theta \ |\ X=x\right]           =  \frac{\displaystyle \int_\Omega \th^{x+1-r}e^{-x\th}dG(\th)}
{\displaystyle \int_\Omega \th^{x-r}e^{-x\th}dG(\th)}. 
\end{equation}

{\bf Example 1}\ Let the prior $G$ be Beta$(v,w)$, $v, w>0$. One can verify 
(see Moll (2015), p.97 for the evaluation of the integrals) that (\ref{bayes}) yields
\begin{equation} \label{bayesB}
\th_{G}(x)   =  
\frac{\displaystyle \sum_{k=0}^{w-1} (-1)^k { w-1 \choose k} \frac{\displaystyle (x-r+v+k)!}{x^{k+1}} \left[ e^x - Exp_{x-r+v+k}(x)\right]}
         {\displaystyle \sum_{k=0}^{w-1} (-1)^k { w-1 \choose k} \frac{\displaystyle (x-r+v+k-1)!}{x^{k}} \left[ e^x - Exp_{x-r+v+k-1}(x)\right]},
 \end{equation}   
where $Exp_j(x):=\sum_{k=0}^{j} x^k/k!$.

{\bf Example 2}\ If the prior $G$ is Uniform$(0,1)$, then (\ref{bayesB}) simplifies to 
\begin{equation} \label{bayesBU}
\th_{G}(x)   
         = 
         \frac{x+1-r}{x} \  \frac{\displaystyle e^x-Exp_{x+1-r}(x)}
         {\displaystyle  e^x-Exp_{x-r}(x)}.
 \end{equation}
We shall adopt the 
empirical Bayes (EB) approach, which relies on the assumption for existence of a prior $G$ which, however, is unknown. Suppose our estimation problem is one in a sequence of similar problems with the same prior distribution.  In this scenario, the results of previous studies can be used to estimate the prior $G$ and/or the Bayes rule $\th_G$ directly.
More precisely, consider a sequence of independent
copies 
\[
(X_1, \Theta_1), (X_2, \Theta_2),\ldots, (X_n, \Theta_n), \ldots 
\]
of the random pair $(X,
\Theta)$, where $\Theta$ has a distribution $G$, and conditional on
$\Theta$, $X$ has the BT distribution (\ref{bt2}). Assume that $X_i$, $i=1,2,\ldots$ are observable, but $\Theta_i$, $i=1,2,\ldots$ are not observable.
We let $X_{n+1}$ stand for the present random observation, and $\underline{X}(n):=(X_1,\ldots,X_n)$ denote the $n$ past observations. Let $\th_{n+1}$ be the present parameter value of the variable $\Theta$. An EB estimator $\th_n(X_{n+1},\underline{X}(n))=: \th_n(X)$ for the parameter $\th$ is a function of the currently observed $X_{n+1}$ and the past data $\underline{X}(n)$.
In general, it is difficult to find an estimator $\theta_n(X)$ for $\theta$ by estimating the Bayes rule $\th_G$ directly. In case of BT distribution, 
Liang (2009) succeeded in  constructing such EB estimator $\th_n(X)$ for $\th_G$ as follows. For $x=r, r+1, \ldots$ let
\[
\psi_n(x):= \frac{1}{n}\sum_{j=1}^n\frac{ c_1(X_j-x)I\{X_j\ge x+1\}}{c_r(X_j)}\quad \mbox{and}\quad q_n(x):=\frac{1}{n}\sum_{j=1}^n \frac{I\{X_j=x\}}{c_r(x)}.
\]
 Define an EB estimator $\th_n(X)$ for each $x=r, r+1, \ldots$ by
\begin{equation} \label{npeb}
\th_n(x):=\min\left\{\frac{\psi_n(x)}{q_n(x)},1\right\}, \qquad q_n(x)\ne 0.
\end{equation}
By definition, the Bayesian estimator $ \th_G(X)$ minimizes the Bayes risk defined (for the squared error loss function) as
\[
R(G,\th_G):=E_{(X,\Theta)}[\Theta-\th_G(X)]^2.
\]
The Bayes risk of the EB estimator $\th_n(X)$ is
\[
R(G,\th_n):=E_nE_{(X_{n+1},\Theta_{n+1})}[\Theta_{n+1}-\th_{n}(X_{n+1})]^2.
\]
The difference 
\[
S(\th_n):=R(G,\th_n)-R(G,\th_G)\ge 0
\]
is called the regret risk of $\th_n$ and measures the quality of $\th_n$. In particular, $\th_n$ is asymptotically optimal for $G$ if $\lim_{n\to \infty}S(\th_n)=0$. Liang (2009) proves that $\th_n$ given by (\ref{npeb}) is asymptotically optimal and studies the rate of convergence to zero of its regret risk $S(\th_n)$ .

\section{Monotone Empirical Bayes Estimator}
 As Van Houwelingen (1977) points out, one issue with the empirical Bayes estimator $\th_n(x)$ is that it is not monotone with respect to $x$  for given values $X_1=x_1,\ldots, X_n=x_n$. 
 On the other hand, it is not difficult to see that the BT distribution (\ref{bt2}) has monotone likelihood ratio (MLR) in $x$, i.e., 
\[
\frac{p_r(x;\theta_2)}{p_r(x;\theta_1)}= \left(\frac{\th_2}{\th_1}\right)^{x-r}e^{-(\th_2-\th_1)}
\]
is an increasing function of $x$ for $0<\th_1<\th_2<1$. Hence, monotonicity is a desirable property for an EB estimator. Estimators for discrete distributions with MLR can be made monotone applying a procedure developed in Van Houwelingen (1977). Consider a simple randomized version of the estimator $\th_n(x)$ represented by the following  function $D(a;x)$ for $a\in [0,1]$:
\[
D(a;x) :=
\left\{
\begin{array}{ll}
        0 & \mbox{if} \ \ \th_n(x) >a, \\
      1 & \mbox{if} \ \ \th_n(x) \le a.
\end{array}
\right.
\]
The number $D(a;x)$ is the probability that an estimate $\th_n(x)$ less than or equal to $a$ is selected if $X=x$.
Hence $D(a,x)$ is a c.d.f. on the action space $(0,1)$ for every $X=x$.
 Define for $a\in [0,1]$ 
\[
\alpha(a) := E (D(a;X))=
\sum_{\{x: \ \th_n(x)\le a \} } p_r(x;a).
\]
Denote 
$F(x;\th) :=  \sum_{k=r}^x p_r(k;\th)$ for $x\ge r$ and $ F(r-1;\th)= 0$.
Now, we can construct a randomized estimator with $D^\ast(a;x)$ as follows
\[
D^\ast(a;x) :=
\left\{
\begin{array}{ll}
        0 & \mbox{if} \ \ \alpha(a)<F(x-1;a) \\
        \frac{\displaystyle \alpha(a)-F(x-1;a)}{\displaystyle F(x;a)-F(x-1;a)} & 
                   \mbox{if} \ \ F(x-1;a)\le \alpha(a)\le F(x;a)\\
          1 & \mbox{if} \ \ F(x;a)<\alpha(a),
\end{array}
\right.
\]
$D^\ast(1;x)=1$, and $D^\ast(0;x)=\lim_{a\downarrow0}D^\ast(a;x)$.  
Let $a\in (\th_0,\th_1)$ be fixed. From the construction of $D^\ast$, it is clear that  $E_aD^\ast(a,X)=\alpha(a)=E_aD(a,X)$.
 It was proven in Van Houwelingen (1977) that $D^\ast$ represents a monotone estimator, which dominates the initial estimator represented by $D$ by having lower Bayes risk, i.e.,  for all $\th\in \Omega$
\[
R(\th,D^\ast)\le R(\th, D).
\] 
Finally, it is not difficult to see that, under the squared error loss function, $D^\ast$ itself is dominated by the non-randomized estimator
\[
\th^\ast_n(x):=\int_0^1 a \, dD^\ast(a;x).
\]
Indeed, using Jensen's inequality, we have
\nbeq
R(\th, \th^\ast_n(X)) & = & E(\th-\th^\ast_n(X))^2 \\
	& = & E \left( \int_0^1 (\th -a)\, dD^\ast (a,X)\right)^2 \\
    & \le & E\left(\int_0^1 (\th -a )^2\, dD^\ast (a,X) \right) \\
    & = & R(\th, D^\ast (a,X)).
    \neeq

\section{Numerical Study}
In practical applications, there is a compelling argument  (Liang (2009)) for $\th$ to take on values in a sub-interval of $(0,1)$. 
Let the prior $\tilde{G}$ be the uniform distribution on $(0.5,0.8)$.
Assuming $r=3$, we find the Bayesian estimator $\th_{\tilde{G}}(X)$ and calculate its 
(minimum) Bayes risk 
\[
E[\th_{\tilde{G}}(X)-\Theta]^2=0.0021.
\]
The maximum likelihood estimator 
$\th_{mle}(X)=(X-3)/X$ has regret risk 
\[
S(\th_{mle}) =E[\th_{mle}(X)-\th_{\tilde{G}}(X)]=0.0935.
\]

\begin{table}[h]
 \centering
\begin{tabular}[c]{ccccc}
\hline 
{\rule{0ex}{3ex} $r$ } & {\rule{0ex}{3ex} $n$} & {\rule{0ex}{3ex}$\hat{S}(\th_n)$ }   &{\rule{0ex}{2ex} $ \hat{S}(\th_n^\ast) $ } &  $S(\th_{mle})$ \\
\hline
3 & {\rule{0ex}{3ex}  100 } & 0.0488  & 0.0242   &  0.0935 \\ 
 &  & (0.0012) &   (0.0008)    &   \\ \hline
3 & {\rule{0ex}{3ex} 500} & 0.0178  & 0.0082  &  0.0935 \\ 
 &  & (0.0004)   & (0.0001)   &   \\ \hline 
\end{tabular}
 \caption{ Estimates for the regret risks of $\th_n$, $\th_n^\ast$, and $\th_{mle}$ (with standard errors in parentheses).}
\end{table}

\begin{figure}[h]
\centering
\includegraphics[width=0.8\textwidth]{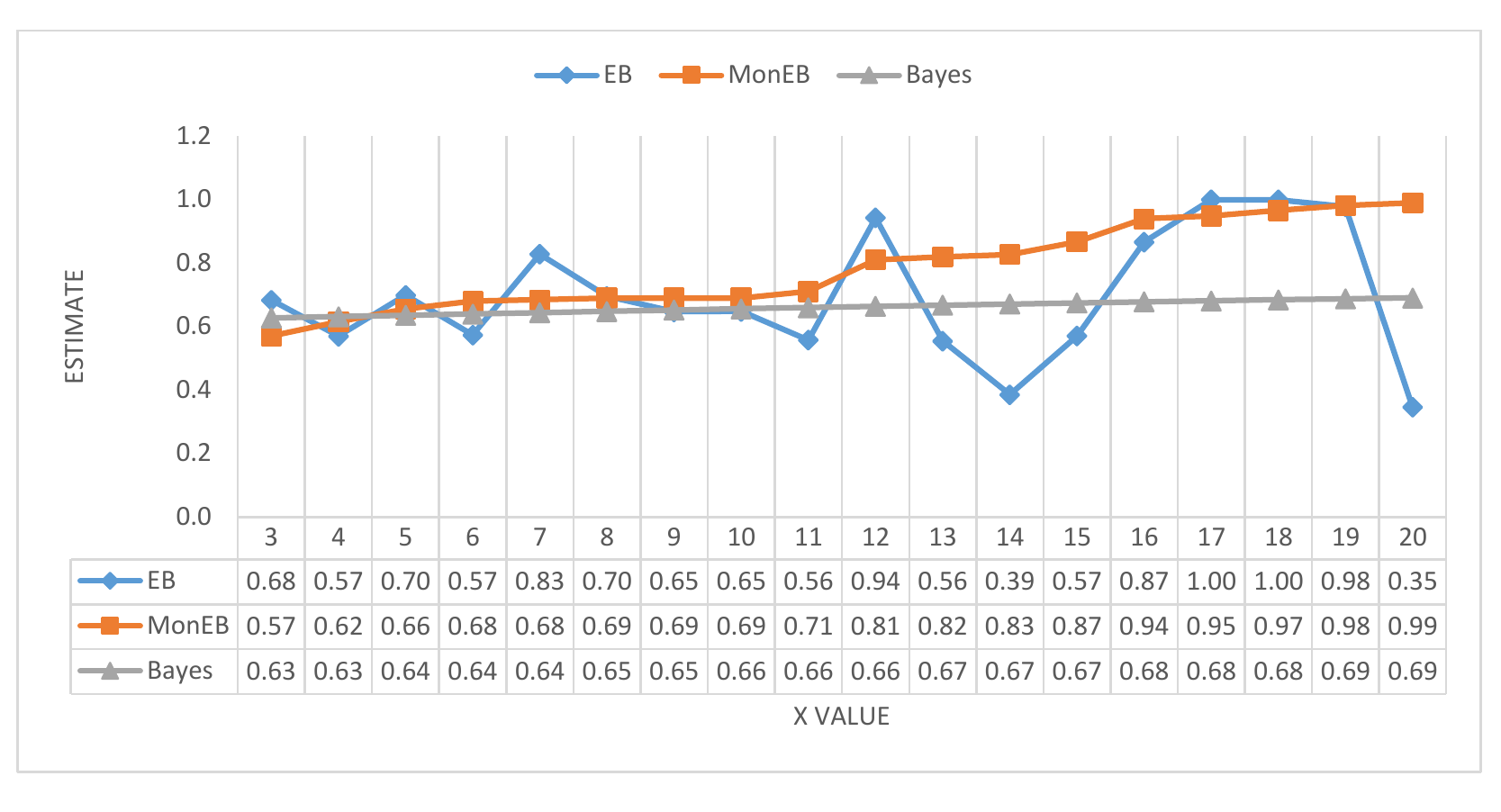}
\caption{\label{fig:1} EB, Monotone EB, and Bayesian estimates based on one simulation for $n=500$, $r=3$, and prior $U(0.5, 0.8)$.}
\end{figure}

\noindent Now adopt the EB framework. 
Consider $n=100$ independent
copies 
\begin{equation} \label{n100}
(X_1, \Theta_1), (X_2, \Theta_2),\ldots ,(X_{100}, \Theta_{100})
\end{equation}
of the random pair $(X, \Theta)$, where $\Theta$ is an uniform $(0.5, 0,8)$ variable and, given $\Theta$, $X$ has the BT distribution (\ref{bt2}). Assume that $X_i$ for  $1\le i\le 100$ are observable, but $\Theta_i$ for  $1\le i\le 100$ are not observable.
For our simulation study, we draw 10 sets like (\ref{n100}). For the  $k^{th}, 1\le k\le 10$, set, the EB estimate $\th^{(k)}_{100}(x)$ is calculated. The value of $S(\th_{100})$ is estimated by the average  (for the 10 samples) $\hat{S}(\th_{100}):={\displaystyle \frac{1}{10}\sum_{k=1}^{10} 
 S(\th^{(k)}_{100})}$ and the standard error is calculated.
 Next, the EB estimator is monotonized and the estimate $\th^{\ast(k)}_{100}(x)$ is computed. 
Similarly to $S(\th_{100})$, we estimate $S(\th_{100}^\ast)$ by the average $\hat{S}(\th_{100}^\ast)$. The entire procedure is 
repeated with $n=500$ in (\ref{n100}). The numerical results are given in Table 1.
The improvement of $\th_n^\ast$ over $\th_n$ is quite substantial. 
It is surprising that even in the case $n=500$, $\th_n$ lacks monotonicity completely. 
To give more insight, the complete results  for one set (\ref{n100}) of size $n=500$ are presented in Figure~1.

\vspace{0.3cm}
{\bf Acknowledgements}\ 
The first author was partially supported by the NFSR at the MES of Bulgaria,
Grant No DFNI-I02/17 while being on leave from the Institute of Mathematics and Informatics at the Bulgarian Academy of Sciences.

\end{document}